\documentclass{amsart}
\usepackage{amssymb}
\usepackage{amsmath}
\usepackage{amsfonts}

\newtheorem{theorem}{Theorem}
\theoremstyle{plain}

\newtheorem{corollary}{Corollary}

\newtheorem{definition}{Definition}

\numberwithin{equation}{section}
\input{tcilatex}

\begin{document}
\title[Modified twisted Dirichlet's type $q$-Euler numbers and polynomials ]{
On $p$-adic interpolating Function Associated with Modified Dirichlet's type
of twisted $q$-Euler Numbers and Polynomials with weight $\alpha $ }
\author{Serkan Araci}
\address{University of Gaziantep, Faculty of Science and Arts, Department of
Mathematics, 27310 Gaziantep, TURKEY}
\email{mtsrkn@hotmail.com}
\author{Mehmet Acikgoz}
\address{University of Gaziantep, Faculty of Science and Arts, Department of
Mathematics, 27310 Gaziantep, TURKEY}
\email{acikgoz@gantep.edu.tr}
\author{Hassan Jolany}
\address{School of Mathematics, Statistics and Computer Science, University
of Tehran, Iran}
\email{hassan.jolany@khayam.ut.ac.ir}
\date{January 18, 2012}

\begin{abstract}
The $q$-calculus theory is a novel theory that is based on finite difference
re-scaling. The rapid development of $q$-calculus has led to the discovery
of new generalizations of $q$-Euler polynomials involving $q$-integers. The
present paper deals with the modified Dirichlet's type of twisted $q$-Euler
polynomials with weight $\alpha $. We apply the method of generating
function and $p$-adic $q$-integral representation on $%
\mathbb{Z}
_{p}$, which are exploited to derive further classes of $q$-Euler numbers
and polynomials. To be more precise we summarize our results as follows, we
obtain some combinatorial relations between modified Dirichlet's type of
twisted $q$-Euler numbers and polynomials with weight $\alpha $. Furthermore
we derive witt's type formula and Distribution formula (Multiplication
theorem) for modified Dirichlet's type of twisted $q$-Euler numbers and
polynomials with weight $\alpha $. In section three, by applying Mellin
transformation we define $q$-analogue of modified twisted $q$-$l$-functions
of Dirichlet's type and also we deduce that it can be written as modified
Dirichlet's type of twisted $q$-Euler polynomials with weight $\alpha $.
Moreover we will find a link between modified twisted Hurwitz-zeta function
and $q$-analogue of modified twisted $q$-$l$-functions of Dirichlet's type
which yields a deeper insight into the effectiveness of this type of
generalizations. In addition we consider $q$-analogue of partial zeta
function and we derive behavior of the modified $q$-Euler $L$-function at $%
s=0$. In final section, we construct $p$-adic twisted Euler $q$-$L$ function
with weight $\alpha $ and interpolate Dirichlet's type of twisted $q$-Euler
polynomials with weight $\alpha $ at negative integers. Our new generating
function possess a number of interesting properties which we state in this
paper.
\end{abstract}

\maketitle

\section{Introduction, Definitions and Notations}

$p$-adic numbers and $L$-functions theory plays a vital and important role
in mathematics. $p$-adic numbers were invented by the German mathematician
Kurt Hensel \cite{hensel}, around the end of the nineteenth century. In
spite of their being already one hundred years old, these numbers are still
today enveloped in an aura of mystery within the scientific community. The $%
p $-adic integral was used in mathematical physics, for instance, the
functional equation of the $q$-Zeta function, $q$-Stirling numbers and $q$%
-Mahler theory of integration with respect to the ring $%
\mathbb{Z}
_{p}$ together with Iwasawa's $p$-adic $q$-$L$ functions. A $p$-adic zeta
function, or more generally a $p$-adic $L$-function, is a function analogous
to the Riemann zeta function, or more general $L$-functions, but whose
domain and target are $p$-adic (where $p$ is a prime number). For example,
the domain could be the $p$-adic integers $%
\mathbb{Z}
_{p}$, a profinite $p$-group, or a $p$-adic family of Galois
representations, and the image could be the $p$-adic numbers $%
\mathbb{Q}
_{p}$ or its algebraic closure. For a prime number $p$ and for a Dirichlet
character defined modulo some integer, the $p$-adic $L$-function was
constructed by interpolating the values of complex analytic $L$-function at
non-positive integers. In this paper our main focus will be on $p$-adic
interpolation of modified Dirichlet's type of twisted $q$-Euler polynomials
with weight $\alpha $. Actually interpolation is the process of defining a
continuous function that takes on specified values at specified points.
During the development of $p$-adic analysis, researches were made to derive
a meromorphic function, defined over the $p$-adic number field, which would
interpolate the same or at least similar values as the Dirichlet $L$%
-function at non-positive integers. Finding the interpolation functions of
special orthogonal numbers and polynomials started by H. Tsumura \cite%
{simsek 2}, and P. T. Young \cite{young}, for the Bernoulli and Euler
polynomials. After Taekyun Kim and Yilmaz Simsek, studied on $p$-adic
interpolation functions of these numbers and polynomials. L. C. Washington 
\cite{Washington}, constructed one-variable $p$-adic-$L$-function which
interpolates generalized classical Bernoulli numbers at negative integers.
Diamond \cite{Diamond}, obtained formulas which express the values of $p$%
-adic $L$-function at positive integers in terms of the $p$-adic log gamma
function. Next Fox in \cite{Fox}, introduced the two-variable $p$-adic $L$%
-functions and T. Kim \cite{Kim 22}, constructed the two-variable $p$-adic $%
q $-$L$-function, which is interpolation function of the generalized $q$%
-Bernoulli polynomials. P. T. Young \cite{young}, gave $p$-adic integral
representations for the two-variable $p$-adic $L$-function introduced by
Fox. T. Kim and S.-H. Rim \cite{Kim 18}, introduced twisted $q$-Euler
numbers and polynomials associated with basic twisted $q$-$l$-functions by
using $p$-adic $q$-invariant integral on $%
\mathbb{Z}
_{p}$ in the fermionic sense. Also, Jang et al. \cite{JAng1}, investigated
the $p$-adic analogue twisted $q$-$l$-function, which interpolates
generalized twisted $q$-Euler numbers attached to Dirichlet's character. In
this paper, we will construct a $p$-adic interpolation function of modified
Dirichlet's type of twisted $q$-Euler polynomials with weight $\alpha $.

Imagine that $p$ be a fixed odd prime number. Throughout this paper $%
\mathbb{Z}
_{p}$, $%
\mathbb{Q}
_{p}$ and $%
\mathbb{C}
_{p}$ will be denote the ring of $p$-adic rational integers, the field of $p$%
-adic rational numbers, and the completion of algebraic closure of $%
\mathbb{Q}
_{p}$, respectively. Let $%
\mathbb{N}
$ be the set of natural numbers and $%
\mathbb{N}
^{\ast }:=%
\mathbb{N}
\cup \left\{ 0\right\} $. In this paper, we assume that $\alpha \in 
\mathbb{Q}
$ and $q\in 
\mathbb{C}
_{p}$ with $\left\vert 1-q\right\vert _{p}<1$.

The $p$-adic absolute value $\left\vert .\right\vert _{p}$, is normally
defined by 
\begin{equation*}
\left\vert x\right\vert _{p}=\frac{1}{p^{r}}\text{,}
\end{equation*}

where $x=p^{r}\frac{s}{t}$ with $\left( p,s\right) =\left( p,t\right)
=\left( s,t\right) =1$ and $r\in 
\mathbb{Q}
$.

As well-known definition, Euler polynomials are defined by%
\begin{equation*}
\frac{2}{e^{t}+1}e^{xt}=\sum_{n=0}^{\infty }\boldsymbol{E}_{n}\left(
x\right) \frac{t^{n}}{n!}=e^{\boldsymbol{E}\left( x\right) t}\text{,}
\end{equation*}

with the usual convention about replacing $\boldsymbol{E}^{n}\left( x\right) 
$ by $\boldsymbol{E}_{n}\left( x\right) $ (for more information, see [8, 9,
14, 15, 16])\bigskip .

A $p$-adic Banach space $B$ is a $%
\mathbb{Q}
_{p}$-vector space with a lattice $B^{0}$ ($%
\mathbb{Z}
_{p}$-module) separated and complete for $p$-adic topology, ie., 
\begin{equation*}
B^{0}\simeq \lim_{\overleftarrow{n\in 
\mathbb{N}
}}B^{0}/p^{n}B^{0}\text{.}
\end{equation*}

For all $x\in B$, there exists $n\in 
\mathbb{Z}
$, such that $x\in p^{n}B^{0}$. Define 
\begin{equation*}
v_{B}\left( x\right) =\sup_{n\in 
\mathbb{N}
\cup \left\{ +\infty \right\} }\left\{ n:x\in p^{n}B^{0}\right\} \text{.}
\end{equation*}

It satisfies the following properties:%
\begin{eqnarray*}
v_{B}\left( x+y\right)  &\geq &\min \left( v_{B}\left( x\right) ,v_{B}\left(
y\right) \right) \text{,} \\
v_{B}\left( \beta x\right)  &=&v_{p}\left( \beta \right) +v_{B}\left(
x\right) \text{, if }\beta \in 
\mathbb{Q}
_{p}\text{.}
\end{eqnarray*}

Then, $\left\Vert x\right\Vert _{B}=p^{-v_{B}\left( x\right) }$ defines a
norm on $B,$ such that $B$ is complete for $\left\Vert .\right\Vert _{B}$
and $B^{0}$ is the unit ball.

A measure on $%
\mathbb{Z}
_{p}$ with values in a $p$-adic Banach space $B$ is a continuous linear map%
\begin{equation*}
f\mapsto \int f\left( x\right) \mu =\int_{%
\mathbb{Z}
_{p}}f\left( x\right) \mu \left( x\right)
\end{equation*}

from $C^{0}\left( 
\mathbb{Z}
_{p},%
\mathbb{C}
_{p}\right) $, (continuous function on $%
\mathbb{Z}
_{p}$) to $B$. We know that the set of locally constant functions from $%
\mathbb{Z}
_{p}$ to $%
\mathbb{Q}
_{p}$ is dense in $C^{0}\left( 
\mathbb{Z}
_{p},%
\mathbb{C}
_{p}\right) $ so.

Explicitly, for all $f\in C^{0}\left( 
\mathbb{Z}
_{p},%
\mathbb{C}
_{p}\right) $, the locally constant functions 
\begin{equation*}
f_{n}=\sum_{i=0}^{p^{n}-1}f\left( i\right) 1_{i+p^{n}%
\mathbb{Z}
_{p}}\rightarrow \text{ }f\text{ in }C^{0}\text{.}
\end{equation*}

Now if ~$\mu \in \boldsymbol{D}_{0}\left( 
\mathbb{Z}
_{p},%
\mathbb{Q}
_{p}\right) $, set $\mu \left( i+p^{n}%
\mathbb{Z}
_{p}\right) =\int_{%
\mathbb{Z}
_{p}}1_{i+p^{n}%
\mathbb{Z}
_{p}}\mu $. Then $\int_{%
\mathbb{Z}
_{p}}f\mu $ is given by the following \textquotedblleft Riemann
sums\textquotedblright 
\begin{equation*}
\int_{%
\mathbb{Z}
_{p}}f\mu =\lim_{n\rightarrow \infty }\sum_{i=0}^{p^{n}-1}f\left( i\right)
\mu \left( i+p^{n}%
\mathbb{Z}
_{p}\right) \text{.}
\end{equation*}

T. Kim defined $\mu $ as follows:%
\begin{equation*}
\mu _{-q}\left( a+p^{n}%
\mathbb{Z}
_{p}\right) =\frac{\left( -q\right) ^{a}}{\left[ p^{n}\right] _{-q}}\text{,}
\end{equation*}

so,%
\begin{eqnarray}
&&I_{-q}\left( f\right)  \label{equation 1} \\
&=&\int_{%
\mathbb{Z}
_{p}}f\left( x\right) d\mu _{-q}\left( x\right)  \notag \\
&=&\lim_{N\rightarrow \infty }\frac{1}{\left[ p^{N}\right] _{-q}}%
\sum_{x=0}^{p^{N}-1}\left( -1\right) ^{x}f\left( x\right) q^{x}\text{, (for
details, see \cite{kim 15},\cite{kim 16},\cite{kim 17}). }  \notag
\end{eqnarray}

Where $\left[ x\right] _{q}$ is a $q$-extension of $x$ which is defined by%
\begin{equation*}
\left[ x\right] _{q}=\frac{1-q^{x}}{1-q}\text{,}
\end{equation*}

note that $\lim_{q\rightarrow 1}\left[ x\right] _{q}=x$ cf. [2-35].

If we take $f_{1}\left( x\right) =f\left( x+1\right) $ in (\ref{equation 1}%
), then we easily see that 
\begin{equation}
qI_{-q}\left( f_{1}\right) +I_{-q}\left( f\right) =\left[ 2\right]
_{q}f\left( 0\right) \text{.}  \label{equation 2}
\end{equation}

By expression (\ref{equation 2}), we readily see that,%
\begin{equation}
\left( -1\right) ^{n-1}I_{-q}\left( f\right) +q^{n}I_{-q}\left( f_{n}\right)
=\left[ 2\right] _{q}\sum_{l=0}^{n-1}\left( -1\right) ^{n-1-l}q^{l}f\left(
l\right) \text{,}  \label{equation 3}
\end{equation}

where $f_{n}(x)=f\left( x+n\right) $.

Recently, Rim et al. \cite{Rim} defined the modified weighted $q$-Euler
numbers $\boldsymbol{E}_{n,q}^{\left( \alpha \right) }$ and the modified
weighted $q$-Euler polynomials $\boldsymbol{E}_{n,q}^{\left( \alpha \right)
}\left( x\right) $ by using $p$-adic $q$-integral on $%
\mathbb{Z}
_{p}$ in the form%
\begin{equation*}
\boldsymbol{E}_{n,q}^{\left( \alpha \right) }=\int_{%
\mathbb{Z}
_{p}}q^{-\xi }\left[ \xi \right] _{q^{\alpha }}^{n}d\mu _{-q}\left( \xi
\right) \text{, for }n\in 
\mathbb{N}
^{\ast }\text{ and }\alpha \in 
\mathbb{Z}
\text{.}
\end{equation*}

Let $C_{p^{n}}=\left\{ w\mid w^{p^{n}}=1\right\} $ be the Cylic group of
order $p^{n}$, and let%
\begin{equation*}
\mathbf{T}_{\mathbf{p}}=\lim_{n\rightarrow \infty }C_{p^{n}}=C_{p^{\infty }}=%
\underset{n\geq 0}{\cup }C_{p^{n}}\text{,}
\end{equation*}%
note that $\mathbf{T}_{\mathbf{p}}$ is locally constant space (for details,
see [17, 23, 24, 27-30, 33]).

In \cite{Park}, let $\chi $ be a Dirichlet's character with conductor $%
d\left( =odd\right) \in 
\mathbb{N}
$ and $w\in \mathbf{T}_{\mathbf{p}}$. If we take $f(x)=\chi \left( x\right)
w^{x}e^{tx}$, then we have $f(x+d)=\chi \left( x\right)
w^{x}w^{d}e^{tx}e^{td}$. $\ $From (\ref{equation 3}), we see that 
\begin{equation}
\int_{X}\chi \left( x\right) w^{x}e^{tx}d\mu _{-q}\left( x\right) =\frac{%
\left[ 2\right] _{q}\sum_{i=0}^{d-1}\left( -1\right) ^{d-1-i}q^{i}\chi
\left( i\right) w^{i}e^{ti}}{q^{d}w^{d}e^{td}+1}\text{.}  \label{equation 4}
\end{equation}

In view of (\ref{equation 4}), it is considered by%
\begin{equation}
F_{w,\chi }^{q}(t)=\frac{\left[ 2\right] _{q}\sum_{i=0}^{d-1}\left(
-1\right) ^{d-1-i}q^{i}\chi \left( i\right) w^{i}e^{ti}}{q^{d}w^{d}e^{td}+1}%
=\sum_{n=0}^{\infty }E_{n,\chi ,w}^{q}\frac{t^{n}}{n!}\text{, }\left\vert
t+\log \left( qw\right) \right\vert <\frac{\pi }{d}\text{.}
\label{equation 5}
\end{equation}

Let us consider the modified twisted $q$-Euler polynomials with weight $%
\alpha $ as follows:%
\begin{equation}
\boldsymbol{E}_{n,q}^{\left( \alpha ,w\right) }\left( x\right) =\int_{%
\mathbb{Z}
_{p}}q^{-\xi }w^{\xi }\left[ x+\xi \right] _{q^{\alpha }}^{n}d\mu
_{-q}\left( \xi \right) \text{, for }n\in 
\mathbb{N}
^{\ast }\text{.}  \label{equation 18}
\end{equation}

By (\ref{equation 18}), and applying combinatorial techniques we have, 
\begin{eqnarray}
\boldsymbol{E}_{n,q}^{\left( \alpha ,w\right) }\left( x\right) 
&=&\sum_{k=0}^{n}\binom{n}{k}q^{\alpha \left( n-k\right) x}\boldsymbol{E}%
_{n-k,q}^{\left( \alpha ,w\right) }\left[ x\right] _{q^{\alpha }}^{k}
\label{equation 19} \\
&=&\sum_{k=0}^{n}\binom{n}{k}q^{\alpha kx}\boldsymbol{E}_{k,q}^{\left(
\alpha ,w\right) }\left[ x\right] _{q^{\alpha }}^{n-k}\text{,}  \notag
\end{eqnarray}

where $\boldsymbol{E}_{n,q}^{\left( \alpha ,w\right) }\left( 0\right) :=%
\boldsymbol{E}_{n,q}^{\left( \alpha ,w\right) }$ are called modified twisted 
$q$-Euler numbers with weight $\alpha $.

By (\ref{equation 18}), we get generating function of modified twisted $q$%
-Euler polynomials as follows:%
\begin{eqnarray}
\tciFourier ^{\left( \alpha \right) }\left( t,x,q,w\right) 
&=&\sum_{n=0}^{\infty }\boldsymbol{E}_{n,q}^{\left( \alpha ,w\right) }\left(
x\right) \frac{t^{n}}{n!}  \label{equation 20} \\
&=&\left[ 2\right] _{q}\sum_{m=0}^{\infty }\left( -1\right) ^{m}w^{m}e^{t%
\left[ x+m\right] _{q^{\alpha }}}\text{.}  \notag
\end{eqnarray}

By using a complex contour integral representation and (\ref{equation 20}),
we get modified twisted Hurwitz-zeta function as follows:%
\begin{eqnarray}
\widetilde{\boldsymbol{\zeta }}_{q}^{\left( \alpha ,w\right) }\left(
s,x\right)  &=&\frac{1}{\Gamma \left( s\right) }\int_{0}^{\infty
}\tciFourier ^{\left( \alpha \right) }\left( -t,x,q,w\right) t^{s-1}dt
\label{equation 21} \\
&=&\left[ 2\right] _{q}\sum_{m=0}^{\infty }\left( -1\right) ^{m}w^{m}\left( 
\frac{1}{\Gamma \left( s\right) }\int_{0}^{\infty }t^{s-1}e^{-t\left[ x+m%
\right] _{q^{\alpha }}}\right)   \notag \\
&=&\left[ 2\right] _{q}\sum_{m=0}^{\infty }\frac{\left( -1\right) ^{m}w^{m}}{%
\left[ m+x\right] _{q^{\alpha }}^{s}}\text{.}  \notag
\end{eqnarray}

By (\ref{equation 20}) and (\ref{equation 21}), we now establish a relation
between $\boldsymbol{E}_{n,q}^{\left( \alpha ,w\right) }\left( x\right) $
and $\widetilde{\boldsymbol{\zeta }}_{q}^{\left( \alpha ,w\right) }\left(
s,x\right) $ as follows:%
\begin{equation}
\widetilde{\boldsymbol{\zeta }}_{q}^{\left( \alpha ,w\right) }\left(
-n,x\right) =\boldsymbol{E}_{n,q}^{\left( \alpha ,w\right) }\left( x\right) 
\text{.}  \label{equation 22}
\end{equation}

In this paper, we construct the generating function of modified Dirichlet's
type twisted $q$-Euler polynomials with weight $\alpha $ in the $p$-adic
case. Also, we give Witt's formula for this type polynomials.\ Moreover, we
obtain a new $p$-adic $q$-Euler $L$-function with weight $\alpha $
associated with Dirichlet's character $\chi ,$ as follows:%
\begin{equation*}
l_{p,q}^{\left( \alpha ,w\right) }\left( -n\mid \chi \right) =\widetilde{E}%
_{n,\chi _{n}}^{\left( \alpha ,w\right) }-\frac{1}{\left[ p^{-1}\right]
_{q^{\alpha F}}^{n}}\chi _{n}\left( p\right) \widetilde{E}_{n,\chi
_{n}}^{\ast \left( \alpha ,w\right) }
\end{equation*}

where $n\in 
\mathbb{N}
^{\ast }$.

\section{\qquad Properties of Modified Dirichlet's type of twisted $q$-Euler
numbers and polynomials}

In this section, by using fermionic $p$-adic $q$-integral equations on $%
\mathbb{Z}
_{p}$, some interesting identities and relations of the modified Dirichlet's
type of twisted $q$-Euler numbers and polynomials with weight $\alpha $, are
given.

\begin{definition}
Let $\chi $ be a Dirichlet's character with conductor $d\left( =odd\right)
\in 
\mathbb{N}
$. For each $n\in 
\mathbb{N}
^{\ast }$ and $w\in T_{p}$. Modified Dirichlet's type of twisted $q$-Euler
polynomials with weight $\alpha $ defined by means of the following
generating function:%
\begin{equation}
\tciFourier ^{\left( \alpha \right) }\left( t,x,q,w\mid \chi \right)
=\sum_{n=0}^{\infty }\widetilde{E}_{n,q}^{\left( \alpha ,w\right) }\left(
x\mid \chi \right) \frac{t^{n}}{n!}  \label{equation 101}
\end{equation}%
where%
\begin{equation}
\tciFourier ^{\left( \alpha \right) }\left( t,x,q,w\mid \chi \right) =\left[
2\right] _{q}\sum_{m=0}^{\infty }\left( -1\right) ^{m}w^{m}\chi \left(
m\right) e^{t\left[ x+m\right] _{q^{\alpha }}}\text{.}  \label{equation 200}
\end{equation}
\end{definition}

From (\ref{equation 101}) and (\ref{equation 200}) we obtain, 
\begin{equation*}
\sum_{n=0}^{\infty }\widetilde{E}_{n,q}^{\left( \alpha ,w\right) }\left(
x\mid \chi \right) \frac{t^{n}}{n!}=\sum_{n=0}^{\infty }\left( \left[ 2%
\right] _{q}\sum_{m=0}^{\infty }\left( -1\right) ^{m}w^{m}\chi \left(
m\right) \left[ x+m\right] _{q^{\alpha }}^{n}\right) \frac{t^{n}}{n!}\text{.}
\end{equation*}

Therefore, we state the following theorem:

\begin{theorem}
Let $\chi $ be a Dirichlet's character with conductor $d\left( =odd\right)
\in 
\mathbb{N}
$. For each $n\in 
\mathbb{N}
^{\ast }$ and $w\in T_{p}$ we have%
\begin{equation}
\widetilde{E}_{n,q}^{\left( \alpha ,w\right) }\left( x\mid \chi \right) =%
\left[ 2\right] _{q}\sum_{m=0}^{\infty }\left( -1\right) ^{m}w^{m}\chi
\left( m\right) \left[ x+m\right] _{q^{\alpha }}^{n}\text{.}
\label{equation 102}
\end{equation}
\end{theorem}

By using (\ref{equation 102}), we can write 
\begin{eqnarray*}
\widetilde{E}_{n,q}^{\left( \alpha ,w\right) }\left( x\mid \chi \right)  &=&%
\left[ 2\right] _{q}\sum_{m=0}^{\infty }\sum_{l=0}^{d-1}\left( -1\right)
^{l+md}w^{l+md}\chi \left( l+md\right) \left[ x+l+md\right] _{q^{\alpha
}}^{n} \\
&=&\frac{\left[ 2\right] _{q}}{\left( 1-q^{\alpha }\right) ^{n}}%
\sum_{l=0}^{d-1}\left( -1\right) ^{l}w^{l}\chi \left( l\right)
\sum_{m=0}^{\infty }\left( -1\right) ^{m}\left( w^{d}\right)
^{m}\sum_{k=0}^{n}\binom{n}{k}\left( -1\right) ^{k}q^{\alpha k\left(
x+l+md\right) } \\
&=&\frac{\left[ 2\right] _{q}}{\left( 1-q^{\alpha }\right) ^{n}}%
\sum_{l=0}^{d-1}\left( -1\right) ^{l}w^{l}\chi \left( l\right) \sum_{k=0}^{n}%
\binom{n}{k}\left( -1\right) ^{k}q^{\alpha k\left( x+l\right)
}\sum_{m=0}^{\infty }\left( -1\right) ^{m}\left( w^{d}\right) ^{m}\left(
q^{\alpha kd}\right) ^{m} \\
&=&\frac{\left[ 2\right] _{q}}{\left( 1-q^{\alpha }\right) ^{n}}%
\sum_{l=0}^{d-1}\left( -1\right) ^{l}w^{l}\chi \left( l\right) \sum_{k=0}^{n}%
\frac{\binom{n}{k}\left( -1\right) ^{k}q^{\alpha k\left( x+l\right) }}{%
q^{\alpha kd}w^{d}+1}\text{.}
\end{eqnarray*}

So, we obtain the following corollary:

\begin{corollary}
Let $\chi $ be a Dirichlet's character with conductor $d\left( =odd\right)
\in 
\mathbb{N}
.$ For each $n\in 
\mathbb{N}
^{\ast }$ and $w\in T_{p}$ we have%
\begin{eqnarray*}
\widetilde{E}_{n,q}^{\left( \alpha ,w\right) }\left( x\mid \chi \right)  &=&%
\left[ 2\right] _{q}\sum_{m=0}^{\infty }\left( -1\right) ^{m}w^{m}\chi
\left( m\right) \left[ x+m\right] _{q^{\alpha }}^{n} \\
&=&\frac{\left[ 2\right] _{q}}{\left[ \alpha \right] _{q}^{n}\left(
1-q\right) ^{n}}\sum_{l=0}^{d-1}\left( -1\right) ^{l}w^{l}\chi \left(
l\right) \sum_{k=0}^{n}\frac{\binom{n}{k}\left( -1\right) ^{k}q^{\alpha
k\left( x+l\right) }}{q^{\alpha kd}w^{d}+1}\text{.}
\end{eqnarray*}
\end{corollary}

By applying $f(\xi )=q^{-\xi }\chi \left( \xi \right) w^{\xi }\left[ x+\xi %
\right] _{q^{\alpha }}^{n}$ into (\ref{equation 1}),%
\begin{eqnarray}
&&\int_{%
\mathbb{Z}
_{p}}q^{-\xi }\chi \left( \xi \right) w^{\xi }\left[ x+\xi \right]
_{q^{\alpha }}^{n}d\mu _{-q}\left( \xi \right)   \label{equation 202} \\
&=&\frac{1}{\left( 1-q^{\alpha }\right) ^{n}}\sum_{k=0}^{n}\binom{n}{k}%
\left( -1\right) ^{k}q^{\alpha kx}\int_{%
\mathbb{Z}
_{p}}\chi \left( \xi \right) w^{\xi }q^{a\xi k-\xi }d\mu _{-q}\left( \xi
\right) \text{,}  \notag
\end{eqnarray}

where from\ (\ref{equation 3}), we easily see that%
\begin{equation}
\int_{%
\mathbb{Z}
_{p}}\chi \left( \xi \right) w^{\xi }q^{\xi \alpha k-\xi }d\mu _{-q}\left(
\xi \right) =\frac{\left[ 2\right] _{q}\sum_{l=0}^{d-1}\left( -1\right)
^{l}q^{\alpha kl}w^{l}\chi \left( l\right) }{q^{\alpha kd}w^{d}+1}\text{.}
\label{equation 203}
\end{equation}

By using (\ref{equation 202}) and (\ref{equation 203}) we obtain%
\begin{eqnarray}
&&\int_{%
\mathbb{Z}
_{p}}q^{-\xi }\chi \left( \xi \right) w^{\xi }\left[ x+\xi \right]
_{q^{\alpha }}^{n}d\mu _{-q}\left( \xi \right)   \notag \\
&=&\frac{1}{\left( 1-q^{\alpha }\right) ^{n}}\sum_{l=0}^{n}\binom{n}{k}%
\left( -1\right) ^{k}q^{\alpha kx}\frac{\left[ 2\right] _{q}\sum_{l=0}^{d-1}%
\left( -1\right) ^{l}q^{\alpha kl}w^{l}\chi \left( l\right) }{q^{\alpha
kd}w^{d}+1}  \notag \\
&=&\widetilde{E}_{n,q}^{\left( \alpha ,w\right) }\left( x\mid \chi \right) 
\text{.}  \label{equation 204}
\end{eqnarray}

Last from equivalent, we obtain Witt's type formula of modified Dirichlet's
type of twisted$\ q$-Euler polynomials with weight $\alpha $ as follows:

\begin{theorem}
Let $\chi $ be a Dirichlet's character with conductor $d\left( =odd\right)
\in 
\mathbb{N}
$. For each $n\in 
\mathbb{N}
^{\ast }$ and $w\in T_{p}$ we obtain%
\begin{equation}
\widetilde{E}_{n,q}^{\left( \alpha ,w\right) }\left( x\mid \chi \right)
=\int_{%
\mathbb{Z}
_{p}}q^{-\xi }\chi \left( \xi \right) w^{\xi }\left[ x+\xi \right]
_{q^{\alpha }}^{n}d\mu _{-q}\left( \xi \right) \text{.}  \label{equation 104}
\end{equation}
\end{theorem}

By (\ref{equation 200}), we obtain $functional$ $equation$ as follows:%
\begin{equation*}
\tciFourier ^{\left( \alpha \right) }\left( t,x,q,w\mid \chi \right) =e^{t
\left[ x\right] _{q^{\alpha }}}\tciFourier ^{\left( \alpha \right) }\left(
q^{x}t,q,w\mid \chi \right) \text{.}
\end{equation*}%
By using the definition of the generating function $\tciFourier ^{\left(
\alpha \right) }\left( t,x,q,w\mid \chi \right) $ as follows:%
\begin{equation*}
\sum_{n=0}^{\infty }\widetilde{E}_{n,q}^{\left( \alpha ,w\right) }\left(
x\mid \chi \right) \frac{t^{n}}{n!}=\left( \sum_{n=0}^{\infty }\left[ x%
\right] _{q^{\alpha }}^{n}\frac{t^{n}}{n!}\right) \left( \sum_{n=0}^{\infty
}q^{n\alpha x}\widetilde{E}_{n,q}^{\left( \alpha ,w\right) }\left( \chi
\right) \frac{t^{n}}{n!}\right) \text{,}
\end{equation*}

by the Cauchy product in the above equation, we have%
\begin{equation*}
\sum_{n=0}^{\infty }\widetilde{E}_{n,q}^{\left( \alpha ,w\right) }\left(
x\mid \chi \right) \frac{t^{n}}{n!}=\sum_{n=0}^{\infty }\left( \sum_{l=0}^{n}%
\binom{n}{l}q^{\alpha lx}\widetilde{E}_{l,q}^{\left( \alpha ,w\right)
}\left( \chi \right) \left[ x\right] _{q^{\alpha }}^{n-l}\right) \frac{t^{n}%
}{n!}\text{,}
\end{equation*}

Therefore, by comparing the coefficients of $\frac{t^{n}}{n!}$ on the both
sides of the above equation, we can state following theorem:

\begin{theorem}
Let $\chi $ be a Dirichlet's character with conductor $d\left( =odd\right)
\in 
\mathbb{N}
$. For each $n\in 
\mathbb{N}
^{\ast }$ and $w\in T_{p}$ we have 
\begin{equation}
\widetilde{E}_{n,q}^{\left( \alpha ,w\right) }\left( x\mid \chi \right)
=\sum_{l=0}^{n}\binom{n}{l}q^{\alpha xl}\widetilde{E}_{l,q}^{\left( \alpha
,w\right) }\left( \chi \right) \left[ x\right] _{q^{\alpha }}^{n-l}\text{.}
\label{equation 10}
\end{equation}
\end{theorem}

So, by using $umbral$ $calculus$ convention in equality (\ref{equation 10}),
we get%
\begin{equation}
\widetilde{E}_{n,q}^{\left( \alpha ,w\right) }\left( x\mid \chi \right)
=\left( q^{\alpha x}\widetilde{E}_{q}^{\left( \alpha ,w\right) }\left( \chi
\right) +\left[ x\right] _{q^{\alpha }}\right) ^{n}\text{,}
\label{equation 11}
\end{equation}

where $\left( \widetilde{E}_{q}^{\left( \alpha ,w\right) }\left( \chi
\right) \right) ^{n}$ is replaced by $\widetilde{E}_{n,q}^{\left( \alpha
,w\right) }\left( \chi \right) $.

From (\ref{equation 3}) we arrive at the following theorem:

\begin{theorem}
Let $\chi $ be a Dirichlet's character with conductor $d\left( =odd\right)
\in 
\mathbb{N}
$, $w\in T_{p}$ and $m\in 
\mathbb{N}
^{\ast }$ we get%
\begin{equation*}
w^{n}\widetilde{E}_{m,q}^{\left( \alpha ,w\right) }\left( n\mid \chi \right)
+\left( -1\right) ^{n-1}\widetilde{E}_{m,q}^{\left( \alpha ,w\right) }\left(
\chi \right) =\left[ 2\right] _{q}\sum_{l=0}^{n-1}\left( -1\right)
^{n-1-l}\chi \left( l\right) w^{l}\left[ l\right] _{q^{\alpha }}^{m}\text{.}
\end{equation*}
\end{theorem}

So, from (\ref{equation 3}), and some combinatorial techniques we can write 
\begin{eqnarray}
\int_{%
\mathbb{Z}
_{p}}q^{-\xi }\chi \left( \xi \right) w^{\xi }\left[ x+\xi \right]
_{q^{\alpha }}^{n}d\mu _{-q}\left( \xi \right)  &=&\frac{\left[ d\right]
_{q^{\alpha }}^{n}}{\left[ d\right] _{-q}}\sum_{a=0}^{d-1}\left( -1\right)
^{a}\chi \left( a\right) w^{a}\int_{%
\mathbb{Z}
_{p}}q^{-d\xi }w^{d\xi }\left[ \frac{x+a}{d}+\xi \right] _{q^{d\alpha
}}^{n}d\mu _{\left( -q\right) ^{d}}\left( \xi \right)   \notag \\
&=&\frac{\left[ d\right] _{q^{\alpha }}^{n}}{\left[ d\right] _{-q}}%
\sum_{a=0}^{d-1}\left( -1\right) ^{a}w^{a}\chi \left( a\right) \boldsymbol{E}%
_{n,q^{d}}^{\left( \alpha ,w^{d}\right) }\left( \frac{x+a}{d}\right) \text{.}
\label{equation 111}
\end{eqnarray}

Therefore, by (\ref{equation 111}), we obtain the following theorem:

\begin{theorem}
Let $\chi $ be a Dirichlet's character with conductor $d\left( =odd\right)
\in 
\mathbb{N}
$, $w\in T_{p}$ and $n\in 
\mathbb{N}
^{\ast }$ we have%
\begin{equation*}
\widetilde{E}_{n,q}^{\left( \alpha ,w\right) }\left( x\mid \chi \right) =%
\frac{\left[ d\right] _{q^{\alpha }}^{n}}{\left[ d\right] _{-q}}%
\sum_{a=0}^{d-1}\left( -1\right) ^{a}w^{a}\chi \left( a\right) \boldsymbol{E}%
_{n,q^{d}}^{\left( \alpha ,w^{d}\right) }\left( \frac{x+a}{d}\right) \text{.}
\end{equation*}
\end{theorem}

\section{Modified Dirichlet's type of twisted $q$-Euler $L$-function with
weight $\protect\alpha $}

In this section, our goal is to consider interpolation function of the
generating functions of modified Dirichlet's type of twisted $q$-Euler
polynomials with weight $\alpha $. For $s\in 
\mathbb{C}
$, $w\in T_{p}$ and $\chi $ be a Dirichlet's character with conductor $%
d(=odd)\in 
\mathbb{N}
$, by applying the Mellin transformation in equation (\ref{equation 200}),
we obtain%
\begin{eqnarray*}
\widetilde{\boldsymbol{L}}_{q}^{\left( \alpha ,w\right) }\left( x,s\mid \chi
\right)  &=&\frac{1}{\Gamma \left( s\right) }\doint t^{s-1}\tciFourier
^{\left( \alpha \right) }\left( -t,x,q,w\mid \chi \right) dt \\
&=&\left[ 2\right] _{q}\sum_{m=0}^{\infty }\left( -1\right) ^{m}w^{m}\chi
\left( m\right) \left( \frac{1}{\Gamma \left( s\right) }\int_{0}^{\infty
}t^{s-1}e^{-t\left[ m+x\right] _{q^{\alpha }}}dt\right) \text{,}
\end{eqnarray*}

so, from above equality, we have%
\begin{equation*}
\widetilde{\boldsymbol{L}}_{q}^{\left( \alpha ,w\right) }\left( x,s\mid \chi
\right) =\left[ 2\right] _{q}\sum_{m=0}^{\infty }\frac{\left( -1\right)
^{m}\chi \left( m\right) w^{m}}{\left[ m+x\right] _{q^{\alpha }}^{s}}\text{.}
\end{equation*}

Consequently, we are in position to define modified Dirichlet's type of
twisted $q$-Euler $L$-function as follows:

\begin{definition}
Let $\chi $ be a Dirichlet's character with conductor $d\left( =odd\right)
\in 
\mathbb{N}
$ and $w\in T_{p}$ we have%
\begin{equation}
\widetilde{\boldsymbol{L}}_{q}^{\left( \alpha ,w\right) }\left( x,s\mid \chi
\right) =\left[ 2\right] _{q}\sum_{m=0}^{\infty }\frac{\left( -1\right)
^{m}\chi \left( m\right) w^{m}}{\left[ m+x\right] _{q^{\alpha }}^{s}}\text{,}
\label{equation 112}
\end{equation}%
for all $s\in 
\mathbb{C}
$. We note that $\widetilde{\boldsymbol{L}}_{q}^{\left( \alpha ,w\right)
}\left( x,s\mid \chi \right) $ is analytic function in the whole complex $s$%
-plane.
\end{definition}

By substituting $s=-n$ into (\ref{equation 112}) we easily get%
\begin{equation*}
\widetilde{\boldsymbol{L}}_{q}^{\left( \alpha ,w\right) }\left( x,-n\mid
\chi \right) =\widetilde{E}_{n,q}^{\left( \alpha ,w\right) }\left( x\mid
\chi \right) \text{,}
\end{equation*}

which led to stating following theorem:

\begin{theorem}
Let $\chi $ be a Dirichlet's character with conductor $d\left( =odd\right)
\in 
\mathbb{N}
$, $w\in T_{p}$ and $n\in 
\mathbb{N}
^{\ast }$, we define 
\begin{equation}
\widetilde{\boldsymbol{L}}_{q}^{\left( \alpha ,w\right) }\left( x,-n\mid
\chi \right) =\widetilde{E}_{n,q}^{\left( \alpha ,w\right) }\left( x\mid
\chi \right) \text{.}  \label{equation 12}
\end{equation}
\end{theorem}

$\widetilde{\boldsymbol{L}}_{q}^{\left( \alpha ,w\right) }\left( 1,s\mid
\chi \right) =\widetilde{\boldsymbol{L}}_{q}^{\left( \alpha ,w\right)
}\left( s\mid \chi \right) $ which is the modified Dirichlet's type of
twisted $q$-Euler $L$-function with weight $\alpha $. Now, by applying
combinatorial techniques we can write,%
\begin{eqnarray}
\widetilde{\boldsymbol{L}}_{q}^{\left( \alpha ,w\right) }\left( s\mid \chi
\right)  &=&\left[ 2\right] _{q}\sum_{m=1}^{\infty }\frac{\left( -1\right)
^{m}\chi \left( m\right) w^{m}}{\left[ m\right] _{q^{\alpha }}^{s}}  \notag
\\
&=&\left[ 2\right] _{q}\sum_{m=1}^{\infty }\sum_{a=0}^{d-1}\frac{\left(
-1\right) ^{a+dm}\chi \left( a+dm\right) w^{a+dm}}{\left[ a+dm\right]
_{q^{\alpha }}^{s}}  \notag \\
&=&\frac{\left[ 2\right] _{q}}{\left[ 2\right] _{q^{d}}}\left[ d\right]
_{q^{\alpha }}^{-s}\sum_{a=0}^{d-1}\left( -1\right) ^{a}\chi \left( a\right)
w^{a}\left[ \left[ 2\right] _{q^{d}}\sum_{m=1}^{\infty }\frac{\left(
-1\right) ^{m}\left( w^{d}\right) ^{m}}{\left[ \left( \frac{a}{d}+m\right) %
\right] _{q^{d\alpha }}^{s}}\right]   \notag \\
&=&\frac{\left[ 2\right] _{q}}{\left[ 2\right] _{q^{d}}}\left[ d\right]
_{q^{\alpha }}^{-s}\sum_{a=0}^{d-1}\left( -1\right) ^{a}\chi \left( a\right)
w^{a}\widetilde{\zeta }_{q^{d}}^{\left( \alpha ,w^{d}\right) }\left( s,\frac{%
a}{d}\right) \text{.}  \label{equation 205}
\end{eqnarray}%
So, by previous calculation we can state following theorem:

\begin{theorem}
Let $\chi $ be a Dirichlet's character with conductor $d\left( =odd\right)
\in 
\mathbb{N}
$ and $w\in T_{p}$ we have%
\begin{equation}
\widetilde{\boldsymbol{L}}_{q}^{\left( \alpha ,w\right) }\left( s\mid \chi
\right) =\frac{\left[ 2\right] _{q}}{\left[ 2\right] _{q^{d}}}\left[ d\right]
_{q^{\alpha }}^{-s}\sum_{a=0}^{d-1}\left( -1\right) ^{a}\chi \left( a\right)
w^{a}\widetilde{\boldsymbol{\zeta }}_{q^{d}}^{\left( \alpha ,w^{d}\right)
}\left( s,\frac{a}{d}\right) \text{.}  \label{equation 206}
\end{equation}
\end{theorem}

We now consider the partial-zeta function $\widetilde{\boldsymbol{H}}%
_{q}^{\left( \alpha \right) }\left( s,a,w\mid F\right) $ as follows:%
\begin{equation}
\widetilde{\boldsymbol{H}}_{q}^{\left( \alpha \right) }\left( s,a,w\mid
F\right) =\left[ 2\right] _{q}\sum_{\underset{m>0}{m\equiv a\left( \func{mod}%
F\right) }}\frac{\left( -1\right) ^{m}w^{m}}{\left[ m\right] _{q^{\alpha
}}^{s}}\text{.}  \label{equation 13}
\end{equation}

If $F\equiv 1\left( \func{mod}2\right) $, then we have 
\begin{eqnarray}
\widetilde{\boldsymbol{H}}_{q}^{\left( \alpha \right) }\left( s,a,w\mid
F\right) &=&\left[ 2\right] _{q}\sum_{\underset{m>0}{m\equiv a\left( \func{%
mod}F\right) }}\frac{\left( -1\right) ^{m}w^{m}}{\left[ m\right] _{q^{\alpha
}}^{s}}  \notag \\
&=&\left[ 2\right] _{q}\sum_{m>0}\frac{\left( -1\right) ^{mF+a}w^{mF+a}}{%
\left[ mF+a\right] _{q^{\alpha }}^{s}}  \notag \\
&=&\frac{\left[ 2\right] _{q}}{\left[ 2\right] _{q^{F}}}\frac{\left(
-1\right) ^{a}w^{a}}{\left[ F\right] _{q^{\alpha }}^{s}}\left[ \left[ 2%
\right] _{q^{F}}\sum_{m>0}\frac{\left( -1\right) ^{m}\left( w^{F}\right) ^{m}%
}{\left[ m+\frac{a}{F}\right] _{q^{\alpha F}}^{s}}\right]  \notag \\
&=&\frac{\left[ 2\right] _{q}}{\left[ 2\right] _{q^{F}}}\frac{\left(
-1\right) ^{a}w^{a}}{\left[ F\right] _{q^{\alpha }}^{s}}\widetilde{%
\boldsymbol{\zeta }}_{q^{F}}^{\left( \alpha ,w^{F}\right) }\left( s,\frac{a}{%
F}\right)  \label{equation 14}
\end{eqnarray}

By expressions (\ref{equation 12}) and (\ref{equation 14}) we get the
following theorem:

\begin{theorem}
Let $F\equiv 1\left( \func{mod}2\right) $, $w\in T_{p}$ , $q$, $s\in 
\mathbb{C}
$, $\left\vert q\right\vert <1$ and $n\in 
\mathbb{N}
^{\ast }$ we have%
\begin{equation}
\widetilde{\boldsymbol{H}}_{q}^{\left( \alpha \right) }\left( -n,a,w\mid
F\right) =\frac{\left[ 2\right] _{q}}{\left[ 2\right] _{q^{F}}}\left(
-1\right) ^{a}w^{a}\left[ F\right] _{q^{\alpha }}^{n}\boldsymbol{E}%
_{n,q^{F}}^{\left( \alpha ,w^{F}\right) }\left( \frac{a}{F}\right) \text{.}
\label{equation 15}
\end{equation}
\end{theorem}

By expressions (\ref{equation 206}) and (\ref{equation 15}), we obtain the
following corollary:

\begin{corollary}
Let $\chi $ be a Dirichlet's character with conductor $d\left( =odd\right)
\in 
\mathbb{N}
$, $w\in T_{p}$ and $F\equiv 1\left( \func{mod}2\right) $ we have%
\begin{equation}
\widetilde{\boldsymbol{L}}_{q}^{\left( \alpha ,w\right) }\left( s\mid \chi
\right) =\sum_{a=0}^{F-1}\chi \left( a\right) \widetilde{\boldsymbol{H}}%
_{q}^{\left( \alpha \right) }\left( s,a,w\mid F\right) \text{.}
\label{equation 23}
\end{equation}
\end{corollary}

By (\ref{equation 19}) and (\ref{equation 15}), we modify the $q$-analogue
of the partial zeta function with weight $\alpha $ as follows:%
\begin{equation}
\widetilde{\boldsymbol{H}}_{q}^{\left( \alpha \right) }\left( s,a,w\mid
F\right) =\frac{\left[ 2\right] _{q}}{\left[ 2\right] _{q^{F}}}\left(
-1\right) ^{a}w^{a}\left[ a\right] _{q^{\alpha }}^{-s}\sum_{l=0}^{\infty }%
\binom{-s}{l}q^{\alpha al}\left( \frac{\left[ F\right] _{q^{\alpha }}}{\left[
a\right] _{q^{\alpha }}}\right) ^{l}\boldsymbol{E}_{l,q^{F}}^{\left( \alpha
,w^{F}\right) }\text{.}  \label{equation 16}
\end{equation}

Let $f\left( \text{=odd}\right) $ and $a$ be the positive integer with $%
0\leq a<F$. Then, (\ref{equation 23}) reduces to%
\begin{equation}
\widetilde{\boldsymbol{L}}_{q}^{\left( \alpha ,w\right) }\left( s\mid \chi
\right) =\frac{\left[ 2\right] _{q}}{\left[ 2\right] _{q^{F}}}%
\sum_{a=0}^{F-1}\chi \left( a\right) \left( -1\right) ^{a}w^{a}\left[ a%
\right] _{q^{\alpha }}^{-s}\sum_{l=0}^{\infty }\binom{-s}{l}q^{\alpha
al}\left( \frac{\left[ F\right] _{q^{\alpha }}}{\left[ a\right] _{q^{\alpha
}}}\right) ^{l}\boldsymbol{E}_{l,q^{F}}^{\left( \alpha ,w^{F}\right) }\text{.%
}  \label{equation 17}
\end{equation}

By expression (\ref{equation 17}), we see that $\widetilde{\boldsymbol{L}}%
_{q}^{\left( \alpha ,w\right) }\left( s\mid \chi \right) $ is an analytic
function $s\in 
\mathbb{C}
$, with except $s=0$. Furthermore, for each $n\in 
\mathbb{Z}
$, with $n\geq 0$, we get%
\begin{equation}
\widetilde{\boldsymbol{L}}_{q}^{\left( \alpha ,w\right) }\left( -n\mid \chi
\right) =\widetilde{E}_{n,q}^{\left( \alpha ,w\right) }\left( \chi \right) 
\text{.}  \label{equation 29}
\end{equation}

By using (\ref{equation 16}), (\ref{equation 17}) and (\ref{equation 29}) we
derive behavior of the modified Dirichlet's type of twisted $q$-Euler $L$%
-function with weight $\alpha $ at $s=0$ as follows:

\begin{theorem}
The following likeable identity%
\begin{equation*}
\widetilde{\boldsymbol{L}}_{q}^{\left( \alpha ,w\right) }\left( 0\mid \chi
\right) =\frac{1+q}{1+w^{F}}\sum_{a=0}^{F-1}\left( -1\right) ^{a}\chi \left(
a\right) w^{a}\text{,}
\end{equation*}%
is true.
\end{theorem}

\section{\textbf{Modified} $\boldsymbol{p}$-\textbf{Adic Twisted} \textbf{%
Interpolation} $q$-$l$-\textbf{Function with weight }$\protect\alpha $}

In this section, we construct modified $p$-adic twisted $q$-Euler $l$%
-function, which interpolate modified Dirichlet's type of twisted $q$-Euler
polynomials at negative integers. Firstly, Washington constructed $p$-adic $%
l $-function which interpolates generalized classical Bernoulli numbers.

Here, we use some the following notations, which will be useful in reminder
of paper.

Let $\omega $ denote the $Kummer$ character by the conductor $f_{\omega }=p$%
. For an arbitrary character $\chi $, we set $\chi _{n}=\chi \omega ^{-n}$, $%
n\in 
\mathbb{Z}
$, in the sense of product of characters.

Let%
\begin{eqnarray*}
\left\langle a\right\rangle  &=&\omega ^{-1}\left( a\right) a=\frac{a}{%
\omega \left( a\right) }\text{,} \\
\left\langle a\right\rangle _{q} &=&\frac{\left[ a\right] _{q}}{\omega
\left( a\right) }\text{.}
\end{eqnarray*}

Thus, we note that $\left\langle a\right\rangle \equiv 1\left( \func{mod}p%
\mathbb{Z}
_{p}\right) $. Let 
\begin{equation*}
A_{j}\left( x\right) =\sum_{n=0}^{\infty }a_{n,j}x^{n}\text{, }a_{n,j}\in 
\mathbb{C}
_{p}\text{,  }j=0,1,2,...
\end{equation*}

be a sequence of power series, each convergent on a fixed subset%
\begin{equation*}
T=\left\{ s\in 
\mathbb{C}
_{p}\mid \left\vert s\right\vert _{p}<p^{-\frac{2-p}{p-1}}\right\} \text{,}
\end{equation*}

of $%
\mathbb{C}
_{p}$ such that

(1) $a_{n,j}\rightarrow a_{n,0}$ as\ $j\rightarrow \infty $ for any $n$;

(2) \ for each $s\in T$ and $\epsilon >0$, there exists an $%
n_{0}=n_{0}\left( s,\epsilon \right) $ such that 
\begin{equation*}
\left\vert \sum_{n\geq n_{0}}a_{n,j}s^{n}\right\vert _{p}<\epsilon \text{
for }\forall j\text{.}
\end{equation*}

So,%
\begin{equation*}
\lim_{j\rightarrow \infty }A_{j}\left( s\right) =A_{0}\left( s\right) \text{%
, for all }s\in T\text{.}
\end{equation*}

This was constructed by Washington \cite{Washington} to indicate that each
functions $\omega ^{-s}\left( a\right) a^{s}$ and 
\begin{equation*}
\sum_{l=0}^{\infty }\binom{s}{l}\left( \frac{F}{a}\right) ^{l}B_{l}\text{,}
\end{equation*}

where $F$ is multiple of $p$ and $f$ and $B_{l}$ is the $l$-th Bernoulli
numbers, is analytic on $T$ (for more information, see \cite{Washington}).

Assume that $\chi $ be a Dirichlet's character with conductor $f\in 
\mathbb{N}
$ with $f\equiv 1(\func{mod}2)$. Thus, we consider the modified Dirichlet's
type of twisted $p$-adic $q$-Euler $l$-function with weight $\alpha $, $%
l_{p,q}^{\left( \alpha ,w\right) }\left( s\mid \chi \right) $, which
interpolate the modified Dirichlet's type of twisted $q$-Euler polynomials
with weight $\alpha $ at negative integers.

For $f\in 
\mathbb{N}
$ with $f\equiv 1\left( \func{mod}2\right) $, let us assume that $F$ is
positive integral multiple of $p$ and $f=f_{\chi }$. We are now ready to
give definition of $l_{p,q}^{\left( \alpha ,w\right) }\left( s\mid \chi
\right) $ as follows:%
\begin{equation}
l_{p,q}^{\left( \alpha ,w\right) }\left( s\mid \chi \right)
=\sum_{a=0}^{F-1}\chi \left( a\right) \left( -1\right) ^{a}w^{a}\left\langle
a\right\rangle _{q^{\alpha }}^{-s}\sum_{l=0}^{\infty }\binom{-s}{l}q^{\alpha
al}\left( \frac{\left[ F\right] _{q^{\alpha }}}{\left[ a\right] _{q^{\alpha
}}}\right) ^{l}\boldsymbol{E}_{l,q^{F}}^{\left( \alpha ,w^{F}\right) }\text{.%
}  \label{equation 24}
\end{equation}

By (\ref{equation 24}), we note that $l_{p,q}^{\left( \alpha ,w\right)
}\left( s\mid \chi \right) $ is analytic for $s\in T$.

For $n\in 
\mathbb{N}
$, we have%
\begin{equation}
\widetilde{E}_{n,\chi _{n}}^{\left( \alpha ,w\right) }=\left[ F\right]
_{q^{\alpha }}^{n}\sum_{a=0}^{F-1}\left( -1\right) ^{a}\chi _{n}\left(
a\right) \boldsymbol{E}_{n,q}^{\left( \alpha ,w\right) }\left( \frac{a}{F}%
\right) \text{.}  \label{equation 25}
\end{equation}

If $\chi _{n}\left( p\right) \neq 0$, then $\left( p,f_{\chi _{n}}\right) =1$%
, and thus the ratio $\frac{F}{p}$ is a multiple of $f_{\chi _{n}}$.

Let%
\begin{equation*}
\lambda =\left\{ \frac{a}{p}\mid a\equiv 0\left( \func{mod}p\right) \text{
for }a_{i}\in 
\mathbb{Z}
\text{ with }0\leq a_{i}<F\right\} \text{.}
\end{equation*}

Thus, we have%
\begin{eqnarray}
&&\left[ F\right] _{q^{\alpha }}^{n}\sum_{\underset{p\mid a}{a=0}%
}^{F-1}\left( -1\right) ^{a}\chi _{n}\left( a\right) \boldsymbol{E}%
_{n,q}^{\left( \alpha ,w\right) }\left( \frac{a}{F}\right) 
\label{equation 26} \\
&=&\frac{1}{\left[ p^{-1}\right] _{q^{\alpha F}}^{n}}\left[ \frac{F}{p}%
\right] _{q^{\alpha }}^{n}\chi _{n}\left( p\right) \sum_{\underset{\eta \in
\lambda }{a=0}}^{\frac{F}{p}}\left( -1\right) ^{\eta }\chi _{n}\left( \eta
\right) \boldsymbol{E}_{n,q}^{\left( \alpha ,w\right) }\left( \frac{\eta }{%
F/p}\right) \text{.}  \notag
\end{eqnarray}

By (\ref{equation 26}), we can define the second modified twisted
generalized Euler numbers attached to $\chi $ as follows:%
\begin{equation}
\widetilde{E}_{n,\chi _{n}}^{\ast \left( \alpha ,w\right) }=\left[ \frac{F}{p%
}\right] _{q^{\alpha }}^{n}\sum_{\underset{\eta \in \lambda }{a=0}}^{\frac{F%
}{p}}\left( -1\right) ^{\eta }\chi _{n}\left( \eta \right) \boldsymbol{E}%
_{n,q}^{\left( \alpha ,w\right) }\left( \frac{\eta }{F/p}\right) \text{.}
\label{equation 27}
\end{equation}

By (\ref{equation 25}), (\ref{equation 26}) and (\ref{equation 27}), we
readily get that%
\begin{eqnarray}
&&\widetilde{E}_{n,\chi _{n}}^{\left( \alpha ,w\right) }-\frac{1}{\left[
p^{-1}\right] _{q^{\alpha F}}^{n}}\chi _{n}\left( p\right) \widetilde{E}%
_{n,\chi _{n}}^{\ast \left( \alpha ,w\right) }  \notag \\
&=&\left[ F\right] _{q^{\alpha }}^{n}\sum_{\underset{p\nshortmid a}{a=0}%
}^{F-1}\left( -1\right) ^{a}\chi _{n}\left( a\right) \boldsymbol{E}%
_{n,q}^{\left( \alpha ,w\right) }\left( \frac{a}{F}\right)
\label{equation 28}
\end{eqnarray}

By (\ref{equation 24}) and (\ref{equation 20}), we readily see that%
\begin{eqnarray*}
&&l_{p,q}^{\left( \alpha ,w\right) }\left( -n\mid \chi \right)  \\
&=&\left[ F\right] _{q^{\alpha }}^{n}\sum_{\underset{p\nshortmid a}{a=0}%
}^{F-1}\left( -1\right) ^{a}\chi _{n}\left( a\right) \boldsymbol{E}%
_{n,q}^{\left( \alpha ,w\right) }\left( \frac{a}{F}\right)  \\
&=&\widetilde{E}_{n,\chi _{n}}^{\left( \alpha ,w\right) }-\frac{1}{\left[
p^{-1}\right] _{q^{\alpha F}}^{n}}\chi _{n}\left( p\right) \widetilde{E}%
_{n,\chi _{n}}^{\ast \left( \alpha ,w\right) }\text{.}
\end{eqnarray*}

Consequently, we state the following Theorem:

\begin{theorem}
Let $n\in 
\mathbb{N}
$, the following equalities%
\begin{equation*}
l_{p,q}^{\left( \alpha ,w\right) }\left( s\mid \chi \right)
=\sum_{a=1}^{F}\chi \left( a\right) \left( -1\right) ^{a}w^{a}\left\langle
a\right\rangle _{q^{\alpha }}^{-s}\sum_{l=0}^{\infty }\binom{-s}{l}q^{\alpha
al}\left( \frac{\left[ F\right] _{q^{\alpha }}}{\left[ a\right] _{q^{\alpha
}}}\right) ^{l}\boldsymbol{E}_{l,q^{F}}^{\left( \alpha ,w^{F}\right) }\text{,%
}
\end{equation*}%
and%
\begin{equation*}
l_{p,q}^{\left( \alpha ,w\right) }\left( -n\mid \chi \right) =\widetilde{E}%
_{n,\chi _{n}}^{\left( \alpha ,w\right) }-\frac{1}{\left[ p^{-1}\right]
_{q^{\alpha F}}^{n}}\chi _{n}\left( p\right) \widetilde{E}_{n,\chi
_{n}}^{\ast \left( \alpha ,w\right) }\text{,}
\end{equation*}%
are true.
\end{theorem}

\end{document}